\documentclass[12pt,twoside]{amsart}
\usepackage{amssymb}
\usepackage{amscd}
\usepackage{xy}

\numberwithin{equation}{section}

\newtheorem{theorem}{Theorem}[section]
\newtheorem{lemma}[theorem]{Lemma}
\newtheorem{proposition}[theorem]{Proposition}
\newtheorem{corollary}[theorem]{Corollary}

\theoremstyle{definition}
\newtheorem{remark}[theorem]{Remark}
\newtheorem{example}[theorem]{Example}

\newcommand{\fm}{{\mathfrak m}}
\newcommand{\fa}{{\mathfrak a}}
\newcommand{\fb}{{\mathfrak b}}

\newcommand{\fq}{{\mathfrak q}}
\newcommand{\fp}{{\mathfrak p}}

\newcommand{\gin}{\operatorname{gin}}
\newcommand{\iin}{\operatorname{in}}
\newcommand{\Hom}{\operatorname{Hom}}
\newcommand{\ER}{\operatorname{Ext}_R}
\newcommand{\Ext}{\operatorname{Ext}}
\newcommand{\HH}{H_{\mathfrak m}}
\newcommand{\hh}{h_{\mathfrak m}}

\newcommand{\depth}{\operatorname{depth}}

\newcommand{\rank}{\operatorname{rank}_K}
\newcommand{\Deg}{\operatorname{Deg}}
\newcommand{\len}{\ell} 

\newcommand{\hdeg}{\operatorname{hdeg}}
\newcommand{\bdeg}{\operatorname{bdeg}}
\newcommand{\reg}{\operatorname{reg}}
 
\newcommand{\moo}{ {\underline {\operatorname{mod}}}\! -\!\! R}

\newcommand{\mif}{\mbox{if} ~}
\newcommand{\s}{\; | \;}
\newcommand{\fall}{\mbox{for all} ~}

\newcommand{\Gm}{\Gamma_{\mathfrak m}} 
\newcommand{\Ass}{\operatorname{Ass}}
\newcommand{\GL}{\operatorname{GL}}
\newcommand{\Ann}{\operatorname{Ann}}

\newcommand {\ZZ}{\mathbb{Z}}

\begin{document}

\title[Comparing Castel\-nuovo-\-Mumford regularity and extended degree]
{Comparing Castel\-nuovo-\-Mumford regularity and extended degree: the
  borderline cases}  

\author[Uwe Nagel]{Uwe Nagel}

\address{Department of Mathematics,
University of Kentucky, 715 Patterson Office Tower, 
Lexington, KY 40506-0027, USA}
\email{uwenagel@ms.uky.edu}

\date{}

\begin{abstract} Castelnuovo-Mumford regularity and any extended
  degree function can be thought of as complexity measures for the
  structure of finitely generated graded modules. A recent result of
  Doering, Gunston, Vasconcelos shows that both can be compared in
  case of a graded algebra. We extend this result to modules and
  analyze when the estimate is in fact an equality. A complete
  classification is 
  obtained if we choose as extended degree the homological or the
  smallest extended degree. The corresponding algebras are
  characterized in three ways: by relations among the algebra
  generators, by using generic initial ideals, and by their Hilbert
  series.  
\end{abstract}


\maketitle
\tableofcontents

\section{Introduction} \label{section-intro} 

In \cite{Vasc-hom-deg} Vasconcelos introduced the concept of an 
extended degree $\Deg (M)$ of a finitely generated graded module $M$
over the polynomial ring $R = K[x_0,\ldots,x_n]$ where $K$ is an
infinite field. It is designed to serve as a  measure of the size and
complexity of the module. Another complexity measure that has attracted a lot of attention is the
Castelnuovo-Mumford regularity $\reg (M)$. In \cite{DG-Vasc}, Theorem
2.4,  it has been shown that both measures can be compared in case of
a standard graded $K$-algebra $A$. In fact, for every extended degree $\Deg$
one has 
\begin{equation} \label{eq-rega}
\reg A \leq \Deg A - 1. 
\end{equation} 
In this note we consider the natural questions how good this bound is and how it 
can  be extended to finitely generated graded $R$-modules.  

Note that a degree shift changes the regularity of the module while its extended degree remains the same. 
 Thus, any comparison between both invariants has to take the degrees
 of the minimal generators of $M$ into account. Denote by $e_R^+(M)$
 the maximum of these  degrees. Then we show as the sought 
 extension  to modules in Theorem \ref{thm-reg-bound} that 
\begin{equation} \label{eq-regm}
\reg M \leq e_R^+(M) + \Deg M - 1. 
\end{equation}
Since a standard graded algebra $A$ is (as $R$-module) generated in
degree zero this specializes to the earlier result (\ref{eq-rega}) in
\cite{DG-Vasc}.  

The key for obtaining estimate (\ref{eq-regm}) is the observation that
for any extended degree the difference $\Deg - \reg$ behaves nicely
with respect to taking general hyperplane sections (cf.\ Theorem
\ref{thm-hypsec}).  

The bulk of this note is devoted to an analysis of when equality holds
in
estimate (\ref{eq-regm}). Note that if $n \geq 2$ 
there are many extended degree functions. Thus, we proceed in two
steps. 
First we derive necessary conditions. This is done for an arbitrary
extended degree in Section \ref{sec-max-reg}. It turns out that  
equality in estimate (\ref{eq-regm}) forces the module to be
cyclic. Thus, we are reduced to analyze equality in estimate 
(\ref{eq-rega}). Theorem \ref{thm-extr-algebras} gives a short list of those 
algebras for which the latter equality has a chance to be true. By
specifying an extended degree this result  allows to obtain  
sufficient conditions for having equality as well. 

This is carried out in two cases.  The first is  Vasconcelos' homological
degree $\hdeg$  
(\cite{Vasc-hom-deg}) which is the prototype of an extended degree.
In  Corollary \ref{cor-hdeg-ex} we characterize the  $K$-algebras
satisfying $\reg (A) = \hdeg (A) - 
1$.  They form 
a proper sublist of the algebras described in  Theorem
\ref{thm-extr-algebras}.  However, we can hope for a more complete converse
of Theorem \ref{thm-extr-algebras} if we choose the extended degree as
small as possible. This is considered in
Section \ref{sec-small}.  
The existence of a smallest extended degree,
denoted by $\bdeg$, has essentially been  established by Gunston 
\cite{Gunst} (cf.\ also Lemma \ref{lem-bdeg-is-small}). It  also has 
the remarkable property that it does not change 
when  passing to the generic initial module (Lemma \ref{lem-gin}). The
main result of this paper characterizes the graded $K$-algebras for
which equality holds in estimate (\ref{eq-rega}) in case $\Deg$ is chosen
to be minimal, i.e. we have 
$\reg (A) = \bdeg (A) - 1$. This  is true if and only if the
algebra is on the list stated in Theorem
\ref{thm-extr-algebras}. Furthermore, we 
show that these algebras can also be characterized  by using generic
initial ideals  or by their 
Hilbert series (cf.\ Theorems \ref{thm-char-s} and
\ref{thm-h-func}). It is rather remarkable that in this case purely
numerical conditions are equivalent to structural properties. 

Another application of Theorem \ref{thm-hypsec} is discussed in
Section \ref{sec-open}. We compare the 
regularity of the modules $M$ and $M/\fq M$ for certain parameter
ideals $\fq$ involving, in particular,  the multiplicity of $\fq$
(cf.\ Proposition \ref{prop-buchs}). This broadly generalizes the main
result of \cite{St-V-mult} for generalized Cohen-Macaulay algebras to
arbitrary modules.  
\smallskip 

\begin{center}
{\bf Acknowledgement}
\end{center}

The author would like to thank W.\ Vasconcelos for helpful comments. 


\section{Comparison with hyperplane section} \label{sec-cdeg} 

The goal of this section is to show a basic inequality about the
difference between any extended degree and the Castelnuovo-Mumford regularity.
Its consequences will be drawn in the following sections. 

We begin by introducing some notation and establishing some preliminary
results.  

Let $R = K[x_0,\ldots,x_n]$ be the polynomial ring over the field $K$ with the standard grading. An {\it extended degree} is a numerical function $\Deg: \moo \to \ZZ$ on the category of finitely generated graded $R$-modules satisfying the following conditions: 
\begin{itemize} 
\item[(i)] If $L = \Gm (M)$ is the submodule of elements of $M$  annihilated by a power of the homogeneous maximal ideal $\fm$ of $R$ then 
$$
\Deg M = \Deg M/L + \len (L) 
$$
where $\len (\_)$ is the length function. 
\item[(ii)] If $l \in R$ is a sufficiently general regular hyperplane
  section of $M$ then  
$$
\Deg M \geq \Deg M/l M. 
$$
\item[(iii)] If $M$ is a Cohen-Macaulay module then 
$$
\Deg M = \deg M 
$$
where $\deg (\_)$ is the ordinary multiplicity function. 
\end{itemize} 

The first example  of an extended degree has been the homological degree $\hdeg$
introduced by Vasconcelos.  
The {\it homological degree} of the module $M$ is
(\cite{Vasc-hom-deg}, Definition 2.8)  
\begin{equation} \label{eq-hom-deg}
\hdeg M = \deg M + \sum_{i=n+2-d}^{n+1} \binom{d-1}{i-n-2+d} \hdeg
(\Ext^i_R(M, R))   
\end{equation} 
where $d = \dim M$. 
Note that the homological degree is defined recursively on the
dimension of the module.  
\smallskip

The Hilbert function of a module $M \in \moo$ is 
$$
h_M: \ZZ \to \ZZ, \; h_M (j) = \rank [M]_j.
$$
It becomes a polynomial in $j$ for $j \gg 0$. In case $d := \dim M >
0$ this Hilbert polynomial has the form  
$$
p_M(j) = \frac{\deg M}{(d-1)!} j^{d-1} + \mbox{lower order terms}. 
$$
If $\dim M = 0$ we set $\deg M = \len (M)$. 

The Hilbert series of $M$ is just the formal power series 
$$
H_M (z) := \sum_{j=0}^{\infty} h_M(j) \cdot z^j. 
$$
\smallskip 

Recall that the irrelevant maximal ideal $(x_0,\ldots,x_n)$ of $R$ is denoted by
$\fm$. The local cohomology functors $\HH^i(\_)$ are the right-derived
functors of $\Gamma_{\fm}(\_)$.  

Let $N$ be any graded $R$-module. Then the end of $N$ is  
$$
e(N) := \sup \{j \in \ZZ \s [N]_j \neq 0\} 
$$  
while its initial degree is 
$$
a (N) = \inf \{j \in \ZZ \s [N]_j \neq 0\}. 
$$
Note that $e(N) = - \infty$ if $N$ is trivial. 
The maximal degree of a
minimal generator of a finitely generated module $N$ is denoted by $e^+_R(N)$, i.e.\  
$$
e^+_R(N) = e (N/ \fm N). 
$$ 
Following \cite{N-compos} we define for $k \geq 0$  integers 
$$
r_k (M) = \max \{i + e(\HH^i(M)) \s i \geq k \}. 
$$
Then $\reg M := r_0(M)$ is the {\it Castelnuovo-Mumford regularity} of
$M$. The number $r_1(M)$ has been called {\it geometric regularity} of
$M$ in \cite{Rossi-TV}. In order to get an upper bound for $r_k(M)$ it
suffices to check finitely many conditions. Indeed, according to
\cite{N-compos}, Theorem 1, we have for $k \geq 1$ that  
$$
r_k(M) \leq m \quad \mbox{if and only if} \quad [\HH^i(M)]_{m+1-i} = 0
\; \fall i \geq k.  
$$
This result can be extended to the case $k = 0$ as follows. 

\begin{lemma} \label{lem-fin-cond} For any integer $m$ we have 
$$
\reg M \leq m \quad \mbox{if and only if} \quad [\HH^i(M)]_{m+1-i} = 0
\; \fall i \geq 0 \; \mbox{and} \; e^+_R (M) \leq m.  
$$
\end{lemma} 

\begin{proof}
Since $e^+_R(M) \leq \reg M$ the conditions are necessary. In order to
show the other implication we induct on the dimension of $M$. If $\dim
M = 0$ we have $M \cong \HH^0 (M)$. Thus, the assumptions imply
$e(\HH^0(M)) \leq m$. Now let $\dim M \geq 1$. Then we conclude as in
the proof of \cite{N-compos}, Theorem 1, for $r_k$ in case $k \geq 1$.  
\end{proof} 

\begin{remark} \label{rem-specializa} 
If $M$  is a graded $K$-algebra $A = R/I$ then the assumption on its
generators can be dropped because $e^+_R(A) = 0 \leq \reg A$. In this
sense,  
Lemma \ref{lem-fin-cond} extends \cite{Rossi-TV}, Corollary 1.2, from
algebras to modules.  
\end{remark} 

Recall that a polynomial $f \in R$ is called {\it $M$-filter regular} if 
$$
f \notin \fp \quad \fall \fp \in \Ass M \setminus \{\fm\}. 
$$
This is equivalent to the condition that the module $0 :_M f$ has finite length.

\begin{remark} \label{rem-seque} 
Let $l \in R$ be a linear $M$-filter regular element. Then it follows that 
$$
\HH^i(M/0 :_M l) \cong \HH^i(M) \quad \fall i \geq 1. 
$$
Hence, the exact sequence induced by multiplication 
$$
0 \to (M/0 :_M l)(-1) \to M \to M/ l M \to 0  
$$
provides the exact sequence 
\begin{equation} \label{eq-seque} 
\begin{CD}
0 \longrightarrow (0 :_M l) (-1) \longrightarrow \HH^0(M) (-1)
\stackrel{l}{\longrightarrow} \HH^0(M)  \longrightarrow \HH^0(M/l M)
\longrightarrow \ldots \\[1ex]  
\ldots \longrightarrow \HH^i(M) \longrightarrow \HH^i (M/l M)
\longrightarrow \HH^{i+1} (M) (-1) \stackrel{l}{\longrightarrow}
\HH^{i+1} (M) \longrightarrow \ldots  
\end{CD} 
\end{equation} 
\end{remark} 
\smallskip

It is well-known that $\reg (M / l M) \leq \reg (M)$ if $l$ is an
$M$-filter regular linear form. Thus, the next 
observation is vacuous if and only if $\reg (M/l M) = \reg (M)$.  

\begin{lemma} \label{lem-h0-gen} 
Let $l \in R$ be an $M$-filter regular linear form. Then we have 
$$
[\HH^0(M)]_j \neq 0 \quad \mif \reg (M/l M) + 1 \leq j \leq \reg (M).
$$
\end{lemma} 

\begin{proof} 
The exact sequence (\ref{eq-seque}) implies for the
geometric regularity  
(cf., for example, \cite{N-compos}, Lemma 2) 
$$
r_1(M) \leq \reg M/l M. 
$$
Hence, if $r := \reg M > \reg M/ l M$ then we must have $[\HH^0(M)]_r
\neq 0$. Now consider the following piece of the sequence (\ref{eq-seque}) 
$$
[\HH^0(M)]_{j-1} \to [\HH^0(M)]_j \to [\HH^0(M/ l M)]_j.
$$
It shows for $j > \reg (M/ lM)$ that $[\HH^0(M)]_j \neq 0$ implies $[\HH^0(M)]_{j-1} \neq 0$. 
Our claim follows. 
\end{proof} 

As a final piece of notation we set  
$$
[\hh^i (M)]_j = \rank [\HH^i (M)]_j. 
$$

Now we are ready for the main result of this section. 

\begin{theorem} \label{thm-hypsec} 
Let $M$ be a finitely generated graded $R$-module of positive
dimension. Suppose $K$ is an infinite field.  
Then we have for every sufficiently general linear form $l \in R$ and any extended degree $\Deg$ 
$$ 
\Deg (M/l M) - \reg (M/l M)  \leq \Deg (M) - \reg (M). 
$$
\end{theorem} 

\begin{proof} 
Since $K$ is infinite, $l$ is $M$-filter regular. 

Let $N$ be any $R$-module. Following \cite{MNP-Bezout} we call $N^{sm}
:= N/\HH^0(N)$ the slight modification of $N$. If it is non-trivial
then it has positive depth.  

For short put $N = M^{sm}$. Comparing the sequences (\ref{eq-seque})
for $M$ and $N$ we obtain the exact sequence  
\begin{equation} \label{eq-H0} 
\begin{CD}
0 @>>> \HH^0(M)/l \HH^0(M) @>>> \HH^0(M/ l M) @>>> \HH^0(N/ l N) @>>> 0. 
\end{CD} 
\end{equation} 
Next, consider the commutative diagram where the vertical maps are
multiplication by $l$  
\begin{equation*}
\begin{CD}
0 @>>> \HH^0 (M) (-1) @>>> M (-1)@>>> N(-1) @>>> 0 \\
& & @VV{l}V @VV{l}V @VV{l}V \\
0 @>>> \HH^0 (M) @>>> M @>>> N @>>> 0.  \\
\end{CD}
\end{equation*} 
The Snake lemma provides the exact sequence 
\begin{equation} \label{eq-M-N} 
\begin{CD}
0 @>>> \HH^0(M)/l \HH^0(M) @>>> M/ l M @>>> N/ l N @>>> 0. 
\end{CD} 
\end{equation} 
Its long exact cohomology sequence induces isomorphisms 
$$
\HH^i(M/ l M) \cong \HH^i(N/l N) \quad \fall i \geq 1.  
$$
Moreover, the sequence (\ref{eq-H0}) implies $e(\HH^0(M/l M)) \geq
e(\HH^0(N/l N))$.  We conclude that  
$$
\reg (N/l N) \leq \reg (M/ l M). 
$$
Comparing the sequences (\ref{eq-H0}) and (\ref{eq-M-N}) we obtain a
graded isomorphism  
$$
(M/l M)^{sm} \cong (N/l N)^{sm}. 
$$
Hence, conditions (i) and (ii) for the extended degree imply 
\begin{eqnarray*} 
\Deg (M/ l M) & = & \len (\HH^0(M/l M)) + \Deg ((M/l M)^{sm}) \\[1ex] 
& = & \len (\HH^0(M/l M)) + \Deg ((N/l N)^{sm}) \\[1ex] 
& \leq & \sum_{j \leq \reg (M/l M)} \left ( [\hh^0(M)]_j + [\hh^0(N/l
  N)]_j \right ) \; + \; \Deg ((N/l N)^{sm}) \\[1ex]  
& = & \sum_{j \leq \reg (M/l M)}  [\hh^0(M)]_j  \; + \; \Deg (N/l N) \\[1ex] 
& \leq & \sum_{j \leq \reg (M/l M)}  [\hh^0(M)]_j  \; + \; \Deg (N) 
\end{eqnarray*} 
where the first estimate is a consequence of sequence (\ref{eq-H0})
and $\reg (N/l N) \leq \reg (M/l M)$ is used to justify the equality
following it. Note that $\depth N > 0$ allows to conclude $\Deg (N/l N) \leq \Deg (N)$ by condition (ii) for the extended degree. 

Now, Lemma \ref{lem-h0-gen} provides the trivial estimate 
$$
\reg (M) - \reg (M/l M) \leq \sum_{j = \reg (M/l M) + 1}^{\reg (M)}
     [\hh^0(M)]_j.  
$$ 
Summing up we obtain 
$$
\reg (M) - \reg (M/l M) + \Deg (M/l M)  \leq \len (\HH^0(M)) + \Deg
(N) = \Deg (M) 
$$ 
which proves the claim. 
\end{proof}

For later use we record when we have equality in Theorem \ref{thm-hypsec} 

\begin{corollary} \label{cor-equal-in-thm} 
Adopting the notation and assumptions of Theorem \ref{thm-hypsec} we have 
$$
\Deg (M/l M) - \reg (M/l M)  = \Deg (M) - \reg (M)  
$$ 
if and only if the following conditions are satisfied \\[-0.6ex]
\begin{itemize} 
\item[(I)] \hspace*{0.5cm} $[\HH^0(M)]_j = [0 :_M l]_j$ \quad for all
  $j < \reg (M/l M)$ \\[-0.6ex] 
\item[(II)] \hspace*{0.5cm} $[\hh^0(M)]_j = 1$ \quad if $\reg (M/l M)
  + 1 \leq j \leq \reg (M)$ \\[-0.6ex] 
\item[(III)] \hspace*{0.5cm} $\Deg (N/ l N) = \Deg (N)$ \quad where $N
  := M/\HH^0(M)$.  
\end{itemize} 
\end{corollary} 

\begin{proof}  
Analyzing the estimates in the proof of Theorem \ref{thm-hypsec} we
see that we have equality in its claim if and only if conditions (II)
and (III) are satisfied and there are equations  
\begin{equation} \label{eq-gleich} 
[\hh^0 (M/l M)]_j = [\hh^0(M)]_j + [\hh^0 (N/l N)]_j \quad \fall j
\leq \reg (M/lM).  
\end{equation}
But the sequences (\ref{eq-seque}) and (\ref{eq-H0}) provide the exact
sequence  
$$
0 \to (0 :_M l)(-1) \to \HH^0(M) (-1) \to \HH^0(M) \to \HH^0(M/ lM)
\to \HH^0 (N/l N) \to 0.   
$$ 
Therefore equations (\ref{eq-gleich}) are equivalent to condition (I). 
\end{proof} 

We now begin to discuss the consequences of Theorem \ref{thm-hypsec}. 

\begin{corollary} \label{cor-Deg-hy} 
If $\Deg (M/l M) = \Deg (M)$ then $\reg (M/l M) = \reg (M)$. 
\end{corollary} 

\begin{proof} 
Using Theorem  \ref{thm-hypsec}, the assumption implies 
$$
\reg (M) \leq \reg (M/ l M). 
$$
But for a general linear form $l$ we always have $\reg (M) \geq \reg
(M/l M)$. 
\end{proof} 

\begin{remark} \label{rem-exdeg-hy} So far, the best known extended
  degree is Vasconcelos' homological degree. An example in \cite{Vasc-hom-deg}, Remark 2.15, shows that there   are modules having positive depth such that  
$$
\hdeg (M/l M) < \hdeg (M). 
$$
Vasconcelos points out that it would be considerably better to have an extended
degree where equality holds in condition (ii). However, it is asking
too much to hope for the existence of an extended degree such that
$\Deg (M/l M) = \Deg (M)$ for general $l$ and every module
$M$. Indeed, there are plenty of modules 
such that  
$$
\reg (M/ lM) < \reg (M). 
$$
\end{remark} 


\section{Modules with maximal regularity} \label{sec-max-reg} 

We will establish a bound on the regularity and then investigate the
modules  whose Castelnuovo-Mumford regularity is maximal with respect
to this bound.  It turns out that this property puts rather severe restrictions on the module structure.  

We will assume that the ground field $K$ is infinite (but confer Remark \ref{rem-finite-field}). 

The following result states the announced regularity bound. 

\begin{theorem} \label{thm-reg-bound} Suppose $K$ is an infinite field. 
Then we have for every finitely generated, graded $R$-module 
$$
\reg (M) \leq e^+_R (M) + \Deg (M) - 1. 
$$
\end{theorem}

\begin{proof}  
Let $l \in R$ be a general linear form. Since $e^+_R (M) = e^+_R (M/ l
M)$, Theorem  \ref{thm-hypsec} shows that, by induction, it suffices to
prove the claim if $\dim M = 0$. But in this case we have $\reg (M) = e
(M)$ and 
\begin{equation} \label{eq-dim-0-case} 
[M]_j \neq 0 \quad \mif \; e^+_R (M) \leq j \leq e (M). 
\end{equation} 
It follows that 
$$
e (M) - e^+_R (M) + 1 \leq \len (M) 
$$ 
and we are done. 
\end{proof} 

If $A = R/I$ is a standard graded $K$-algebra then $e^+_R(A) =
0$. Therefore, the theorem above has the following consequence which
has first been shown as  Theorem
2.4 in \cite{DG-Vasc}.   

\begin{corollary} \label{cor-reg-b} 
Let $A$ be a standard graded $K$-algebra where $K$ is an infinite
field. Then 
$$
\reg (A) \leq \Deg (A) - 1. 
$$
\end{corollary}  

The goal of the remainder of this section is to study  the
modules which have maximal regularity with respect to Theorem
\ref{thm-reg-bound}.  

The next observation shows that it suffices to characterize equality
for algebras. 

\begin{lemma} \label{lem-extr-is-cyclic} 
Let $M$ be a  finitely generated, graded $R$-module such that 
$$
\reg (M) = e^+_R (M) + \Deg (M) - 1. 
$$ 
Suppose $K$ is infinite and 
let $l \in R$ be a general linear form. Then we have: 
\begin{itemize} 
\item[(a)] If $\dim M \geq 1$ then 
$$
\reg (M/l M) = e^+_R (M/l M) + \Deg (M/l M) - 1.
$$
\item[(b)] $M$ is a cyclic $R$-module. 
\end{itemize} 
\end{lemma}  

\begin{proof} 
We will use twice the fact that the degrees of the minimal generators of $M$
and of $M/l M$ coincide for general $l$. First, it shows that claim (a) is a
consequence of Theorem \ref{thm-reg-bound} and Theorem \ref{thm-hypsec}. 

Second, using induction on $\dim M$ and part (a) we see that it
suffices to show claim (b) if $\dim (M) = 0$. In this case the assumption
becomes 
$$
e (M) - e^+_R (M) + 1 = \len (M). 
$$
Hence condition (\ref{eq-dim-0-case}) implies  
\begin{equation} \label{eq-extr-dim-0}
h_M (j) = 1 \quad \mif \; a(M) = e^+_R(M) \leq j \leq e(M) 
\end{equation} 
which in particular proves the claim. 
\end{proof} 

We need some preparation for the next result. The graded $K$ dual of
an $R$-module $N$ is denoted by $M^{\vee} = \oplus_{j \in \ZZ} \Hom_K
([M]_{-j}, K)$. By local duality  
$$
\HH^i (M)^{\vee} \cong \ER^{n+1-i} (M, R)(-n-1) 
$$
is a finitely generated graded $R$-module. 

\begin{lemma} \label{lem-ex-ideal} In $R = K[x_0,\ldots,x_n]$ consider
  the ideal  
\begin{equation*} \label{eq-ideal} 
I = (f_n l_0, f_n f_{n-1} l_1,\ldots, f_n \ldots f_{t+1} l_{n-t-1},
f_n \ldots f_t) 
\end{equation*} 
where $0 \leq t \leq n$, every $f_i \neq 0$ is a homogenous polynomial
of degree $d_i \geq 0,\; d_n, d_t \geq 1$, and  every $l_i$ is a
linear form.  (Note that in case $n = t$ the  ideal
$I$ is simply defined as $I = (f_n)$.) 
If $I$ has (as indicated) $n+1-t$ minimal generators then the local
cohomology modules of $A = R/I$ are  
$$
\HH^i(A)^{\vee} \cong \left \{ 
\begin{array}{ll}
(R/(l_0,\ldots,l_{n-1-i}, f_i))(d_i + \ldots + d_n - i - 1) & \mif \;
  t \leq i \leq n \\ 
0 & \mbox{otherwise}. 
\end{array} \right.
$$ 
In particular, $\HH^i(A)^{\vee}$ either vanishes or is (up to a degree
shift) a hypersurface ring of dimension $i$.  

Furthermore, the Hilbert series of $A$ is 
$$
H_A (z) = \sum_{j=t}^n \frac{1-z^{d_j}}{(1-z)^{j+1}} \prod_{i=j+1}^n z^{d_i} 
$$ 
where the product is defined to be $1$ if $j+1 > n$. 
\end{lemma} 

\begin{proof} 
Since $I$ has $n+1-t$ minimal generators we see that 
\begin{eqnarray*}
l_j & \notin & (l_0,\ldots,l_{j-1}) \quad \mif \; 0 \leq j \leq t \\[.6ex]
f_j & \notin & (l_0,\ldots,l_{n-j-1}) \quad \mif \; t \leq j \leq n. 
\end{eqnarray*} 
It follows that the rings $R/(l_0,\ldots,l_{n-1-i}, f_i)$ have the
asserted dimension $i$ and that $f_n$ is the greatest common divisor
of the polynomials in $I$. Thus, we can write $I = (f_n) \cap \fb$
where $\fb \subset R$ is an ideal having codimension at least two or
$\fb = R$. Consider the exact sequence  
$$
0 \to (f_n)/\fb \cap (f_n) \to R/I \to R/(f_n) \to 0. 
$$
If $I$ is not a principal ideal then the dimension of $(f_n)/\fb \cap
(f_n) \cong (\fb + (f_n))/\fb$ is less than $\dim R/I$. Thus, using
the long exact cohomology sequence we get in any case 
$$
\HH^n(A) \cong \HH^n(R/ f_n R) \cong (R/f_n R)^{\vee} (-d_n + n + 1)
$$ 
as claimed. (The first isomorphism is also a consequence of the more
general Proposition 3.7 in \cite{BN}.) 

In order to compute the cohomology of $A$ we induct on $n$. If $n = 0$
then $I = (f_n)$ and the last computation gives the claim.  
Let $n \geq 1$. Due to the last isomorphism it suffices to compute the
$i$-th cohomology module for $i \leq n-1$.  

Multiplication by $f_n$ on $A$ induces the exact sequence 
\begin{equation} \label{eq-seq}
0 \to (R/I : f_n)(-d_n) \to R/I \to R/f_n R \to 0. 
\end{equation}
Since $R/f_n R$ is Cohen-Macaulay of dimension $n$ it provides isomorphisms 
\begin{equation} \label{eq-iso} 
\HH^i(A) \cong \HH^i(R/I : f_n)(-d_n) \quad (i \leq n-1). 
\end{equation} 
Consider the following ideal  in $R$ 
$$
J = (f_{n-1} l_1,\ldots, f_{n-1} \ldots f_{t+1} l_{n-t-1}, f_{n-1} \ldots f_t). 
$$ 
We have 
$$
I : f_n = l_0 R + J. 
$$
Thus,  putting $\bar{R} = R/l_0 R$ we get the isomorphism
\begin{equation} \label{eq-isored}
R/I : f_n \cong \bar{R}/J \bar{R}. 
\end{equation}
But the induction hypothesis applies to $J \bar{R}$. Hence, the isomorphisms (\ref{eq-iso}) provide the claim about the cohomology of $A$. 

It remains to determine the Hilbert series of $A$. If $n = 0$ we obtain 
$$
H_A (z) = \sum_{j=0}^{d_0 - 1} z^j = \frac{1-z^{d_0}}{1-z}
$$ 
as claimed. Let $n \geq 1$. If $t = n$ then $A$ is a hypersurface ring of degree $d_n$ and its Hilbert series is $H_A (z) =  \frac{1-z^{d_n}}{(1-z)^{n+1}}$. Let $t < n$.  
Then the sequence (\ref{eq-seq}) and the isomorphism (\ref{eq-isored}) together with the induction hypothesis provide 
\begin{eqnarray*}
H_A (z) & = & \frac{1-z^{d_n}}{(1-z)^{n+1}} + z^{d_n} \left [ \sum_{j=t}^{n-1}  \frac{1-z^{d_j}}{(1-z)^{j+1}} \prod_{i=j+1}^{n-1} z^{d_i} \right ] \\
 & = & \sum_{j=t}^n \frac{1-z^{d_j}}{(1-z)^{j+1}} \prod_{i=j+1}^n z^{d_i} 
\end{eqnarray*} 
completing the proof. 
\end{proof} 

We are ready for the main result of this section. 

\begin{theorem} \label{thm-extr-algebras} 
Let $A \neq 0$ be a standard graded $K$-algebra where $K$ is an infinite
field. If 
$$ 
\reg (A) = \Deg (A) - 1
$$ 
for some extended degree $\Deg$ then $A$ is (as $K$-algebra) isomorphic to $R/I$ where $R = K[x_0,\ldots,x_n]$, $n = \dim A$, and 
$$
I = (f_n x_0, f_n f_{n-1} x_1,\ldots, f_n \ldots f_{t+1} x_{n-t-1}, f_n \ldots f_t)
$$
with homogenous polynomials  $f_i \neq 0$  of degree $d_i \geq 0$, $t = \depth A$,  and $d_n, d_t \geq 1$ such that $f_i \notin (x_0,\ldots,x_{n-i-1})$ for all $i=t,\ldots,n-1$.
\end{theorem} 

\begin{proof} We
 will proceed in several steps. First, we show several properties of
 the algebras in question,  second, we use induction  to  conclude the
 proof. 
\smallskip 

{\it Step 1}. Assume $\dim A = 0$. Then we will show that $A \cong K[x]/x^d$ for
some integer $d \geq 0$. 

Indeed, the equations (\ref{eq-extr-dim-0}) read in this case as 
$$
h_A (j) = 1 \quad \mif \; 0 \leq j \leq e(A). 
$$ 
In particular, $\rank [A]_1 \leq 1$ implies that $A$ is isomorphic to
a quotient of a polynomial ring in one variable. 
The claim of Step 1 follows. 
\smallskip 

{\it Step 2}. Put $B = A/\HH^0 (A)$. Then we claim that  
\begin{equation} \label{eq-st-2}
\reg (B) \leq \reg (A) - \len (\HH^0 (A)). 
\end{equation}
This assertion is trivially true if $\depth A > 0$. Thus, assume $\HH^0 (A)
\neq 0$. Then the assumption and Theorem \ref{thm-reg-bound} provide 
\begin{equation} 
\reg (A) + 1 = \Deg (A) = \Deg (B) + \len (\HH^0 (A)) \geq \reg (B) +
1 +  \len (\HH^0 (A)) 
\end{equation} 
and inequality (\ref{eq-st-2}) is shown. 
\smallskip 

{\it Step 3}. Suppose $\dim A \geq 1$ and $\HH^0 (A) \neq 0$. Then we
show that 
\begin{equation} \label{eq-st-3-1}
\len (\HH^0 (A)) \leq e(\HH^0 (A)) - a(\HH^0 (A)) + 1. 
\end{equation}
In order to do that write $A = R/I$ and denote the saturation of $I$
by $I^{sat}$. There are isomorphisms $\HH^0 (A) \cong I^{sat}/I$ and $B \cong
R/I^{sat}$. Estimate (\ref{eq-st-2}) provides $\reg (A) > \reg (B) =
r_1 (A)$ which 
implies $\reg (A) = e(\HH^0 (A))$. 
Since $\reg (I^{sat}) = \reg (B) + 1$ 
it follows  using estimate (\ref{eq-st-2}) again 
\begin{equation} \label{eq-st-3-2}
\left \{ \quad \begin{split}
a (\HH^0 (A)) \leq e^+_R(I^{sat}) & \leq  \reg (B) + 1 \\
& \leq  \reg (A) + 1 - \len (\HH^0 (A)) \\
& =  e(\HH^0 (A)) + 1 - \len (\HH^0 (A)) 
\end{split} \right.
\end{equation}
which gives our claimed estimate (\ref{eq-st-3-1}). 
\smallskip 

{\it Step 4}. If $\HH^0 (A)$ is non-trivial then we claim that  it is
a cyclic $R$-module and  
$$
[h_{\fm}^0 (A)]_j = 1 \quad \mif \; a(\HH^0 (A)) \leq j \leq e(\HH^0 (A)). 
$$
We argue by induction on $\dim A$. If $\dim A = 0$ this claim is shown
in Step 1. Let $\dim A > 0$ and let $l$ be a general linear form. Then
the sequence (\ref{eq-seque}) provides an embedding  
$$
\HH^0(A)/ l \HH^0(A) \hookrightarrow \HH^0 (A/l A). 
$$
According to Lemma \ref{lem-extr-is-cyclic}, (a), the induction
hypothesis applies to $A/l A$. Hence \\ 
 $\HH^0(A)/ l \HH^0(A)$ and thus $\HH^0 (A)$ must be a cyclic
module. This implies  
$$
 e(\HH^0 (A)) - a(\HH^0 (A)) + 1 \leq \len (\HH^0 (A)). 
$$ 
Comparing with estimate (\ref{eq-st-3-1}) we get
\begin{equation} \label{eq-st-4} 
e(\HH^0 (A)) - a(\HH^0 (A)) + 1  =  \len (\HH^0 (A))
\end{equation}
 which implies the claim of this step about the Hilbert function of $\HH^0(A)$.  
\smallskip 

{\it Step 5}. Now we can show 
\begin{itemize} 
\item[(i)] If $\HH^0 (A) \neq 0$ then $a (\HH^0(A)) = e^+_R (I^{sat})
  = r_1 (A) + 1$.  
\item[(ii)] If $\dim A > 0$ then $\reg (B) = \Deg (B) - 1$. 
\end{itemize} 
Indeed,  equation (\ref{eq-st-4})  
implies  equality for all the estimates in (\ref{eq-st-3-2}).  This shows
claim (i) and  
$$
\reg (B) = \reg (A) - \len (\HH^0 (A)). 
$$
Since by assumption 
$$
\reg (A) = \Deg (A) - 1 = \Deg (B) + \len (\HH^0(A)) - 1 
$$
claim (ii) follows. 
\smallskip 

{\it Step 6}. Suppose temporarily that $\HH^0 (A) \cong I^{sat}/I$ is
non-trivial. Then claim (i) in Step 5 says that there is a minimal
generator of $I^{sat}$, say $f$, of maximal degree whose residue class
generates $\HH^0 (A)$ as $R$-module. Combined with Step 4 we conclude
that in any case there is a homogeneous ideal $\fa$ such that  
\begin{equation} 
I^{sat} = \fa + f R \quad \mbox{and} \quad I = \fa + f (l_0^k, l_1,\ldots,l_n) 
\end{equation} 
for some integer $k \geq 0$ where $\{l_0, l_1,\ldots,l_n\}$ is a
regular sequence of linear forms.  (Note that $I = I^{sat}$ if and
only if $k = 0$.)  
\smallskip 
 
{\it Step 7}. Now we are ready to establish the assertion of Theorem
\ref{thm-extr-algebras}. We will 
use Step 6 and argue by induction on $\dim A$. If $\dim A = 0$ the
claim is shown in Step 1. Let $\dim A > 0$. Write $A = R/I$.  
First, suppose that $\depth A > 0$. Let $l$ be a general linear
form. After a change of coordinates we may assume $l = x_n$.  
Lemma \ref{lem-extr-is-cyclic} shows that the induction hypothesis
applies to $A/l A \cong \bar{R}/\bar{I}$ where $\bar{R} = R/x_n R$.
Since $l$ is not a zero divisor for $A$ the minimal generators of
$\bar{R}$ lift to minimal generators of $I$ preserving their
factorization. The claim of our statement follows. (Note that $\depth
A/l A = \depth A - 1$.)  

Second, suppose $\depth A = 0$. Then $t := \depth R/I^{sat} \geq
1$. Hence, by Step 5, (ii) and the first part of Step 7 we may assume that
$I^{sat}$ is of the form as claimed. In particular, $f_n \cdots f_t$
is the unique minimal generator of largest degree of $I^{sat}$ if
$\deg f_t \geq 2$. If $\deg f_t = 1$ we may still assume that we have
in Step 6 $f = f_n \ldots f_t$ and  
$$
\fa = (f_n x_0, f_n f_{n-1} x_1,\ldots, f_n \ldots f_{t+1} x_{n-t-1}) 
$$
such that $I = \fa + f (f_0, l_1,\ldots,l_n)$ with suitable linear forms $l_i$ and a form $f_0 \neq 0$ of degree $d_0 \geq 1$. Since $f (x_0,\ldots,x_{n-1-t}) \subset \fa$ we get 
$$
(x_0,\ldots,x_{n-1-t}) \subset \Ann \HH^0(A) = (f_0, l_1,\ldots,l_n). 
$$ 
Thus, possibly after a change of coordinates we may assume 
$$
(f_0, l_1,\ldots,l_n) = (x_0,\ldots,x_{n-1}, f_0). 
$$
Then $I = \fa + f (x_0,\ldots,x_{n-1}, f_0)$ is of the form as claimed. 
\end{proof} 

Combined with Lemma \ref{lem-extr-is-cyclic} we obtain 

\begin{corollary} \label{cor-ex-mod} 
Suppose the ground field $K$ is infinite. Let $M \neq 0$ be a finitely generated graded $R$-module such that 
$$
\reg (M) = e^+_R (M) + \Deg (M) - 1 
$$ 
for some extended degree $\Deg$. 
Then $M$ is a cyclic $R$-module which is up to a degree shift as graded $K$-algebra isomorphic to one of the algebras in Theorem \ref{thm-extr-algebras}. 
\end{corollary} 

Using as extended degree Vasconcelos' homological degree (cf.\ (\ref{eq-hom-deg})) we get a complete characterization. 

\begin{corollary} \label{cor-hdeg-ex} 
Let $K$ be an infinite field and let  $M \neq 0$ be a finitely generated graded $R$-module such that  
$$
\reg (M) = e^+_R (M) + \hdeg (M) - 1.  
$$ 
Then $M$ is a cyclic $R$-module which is up to a degree shift as graded $K$-algebra isomorphic to $R/I$ where 
$$
I = (f_n x_0, f_n f_{n-1} (x_1,\ldots,x_{n-1}, f_0))
$$
with homogeneous polynomials $f_i \neq 0$ of degree $d_i \geq 0$, $d_n \geq 1$, $f_{n-1} \notin (x_0)$, and $f_0 \notin (x_0,\ldots,x_{n-1})$. 
\end{corollary} 

\begin{proof} 
According to the previous result it suffices to compare the regularity and the homological degree of the algebras in Theorem \ref{thm-extr-algebras}. Let $A$ be such an algebra. Put $t = \depth A$. Note that $\dim A = n$. Then Lemma \ref{lem-ex-ideal} provides 
$$
\reg (A) = d_t + \ldots + d_n  - 1 
$$
and 
$$
\hdeg (A) = d_n + \sum_{i=t}^{n-1} \binom{n-1}{i} d_i. 
$$
Therefore, we get $\reg (A) = \hdeg (A) - 1$ if and only if $d_i = 0$ for all $i$ such that $1 \leq i \leq n-2$ proving the claim. 
\end{proof} 

In the next section we will show that there is an extended degree
function $\bdeg$ such that we have for {\em every} algebra $A$ as in
Theorem \ref{thm-extr-algebras}   
$$
\reg (A) = \bdeg (A) - 1 \; ( = d_t + \ldots + d_n  - 1). 
$$ 


\section{The smallest extended degree} \label{sec-small} 

In his thesis \cite{Gunst} T.\ Gunston introduced a very particular
extended degree which can be described as the smallest among all
extended degrees. Using this degree we obtain several
characterizations of the algebras occurring in Theorem
\ref{thm-extr-algebras}.  

Let us begin with the axiomatic description of this degree (cf.\ \cite{Gunst}, Theorem 3.1.2). 

\begin{theorem} \label{thm-bdeg} 
Let $K$ be an infinite field. Then there is a unique numerical function $\bdeg: \moo \to \ZZ$ satisfying the following conditions: 
\begin{itemize} 
\item[(1)] If $L = \Gamma_{\fm} (M)$ then 
$$
\bdeg (M) = \bdeg (M/L) + \len (L). 
$$
\item[(2)] If $l \in R$ is a sufficiently general regular hyperplane section of $M$ then 
$$
\bdeg M = \bdeg M/l M. 
$$
\item[(3)] $\bdeg (0) = 0$.
\end{itemize} 
\end{theorem} 

T.\ Gunston called the function $\bdeg$ the extremal cohomological
degree. It is easy to see that it is an extended degree and in this
note we prefer to call the function the {\it smallest extended
  degree}. This is justified because of the next observation.  

\begin{lemma} \label{lem-bdeg-is-small} 
If the field $K$ is infinite then we have for every extended degree
$\Deg$ on $\moo$  
$$
\bdeg (M) \leq \Deg (M) \quad \fall M \in \moo. 
$$
\end{lemma} 

\begin{proof} 
The result is contained in T.\ Gunston's thesis \cite{Gunst}. However,
for the convenience of the reader we include the short argument using
the axiomatic description only.  

We induct on $\dim M$. If $\dim M = 0$ then $M$ has finite length and 
$$ 
\bdeg (M) = \Deg (M) = \len (M). 
$$
Let $\dim M > 0$. Then $N := M/\HH^0 (M)$ has positive depth and by induction we have for a general linear form $l \in R$ that $\bdeg (N/lN) \leq \Deg (N/l N)$. Thus, it follows 
\begin{eqnarray*} 
\bdeg (M) = \bdeg (N) + \len(\HH^0 (M)) & = & \bdeg (N/l N) + \len (\HH^0(M)) \\[1ex]
& \leq & \Deg (N/l N) + \len(\HH^0(M)) \\[1ex]
& \leq & \Deg (N) + \len(\HH^0(M)) = \Deg (M). 
\end{eqnarray*}
\end{proof} 

Now we are touching upon the theory of Gr\"obner bases where we assume that the field $K$ is infinite. Let $F = \oplus_{i=0}^r R e_i$ be a free graded $R$-module of rank $r$. We will always use the {\em reverse lexicographic order} on the monomials of $F$, i.e.\ if $m, n \in R$ are monomials then $m e_i > n e_j$ if either $\deg m e_i > \deg n e_j$ or the degrees are the same and $m > n$ in the reverse lexicographic order of the monomials in $R$ or $m = n$ and $i < j$. 

Denote by $\GL(n+1)$ the group of $K$-linear graded automorphisms of $R$ and let $\GL(F)$ be the group of $R$-linear graded automorphisms of $F$. Then $G := \GL(n+1) \ltimes \GL(F)$ acts on $F$ through $K$-linear graded automorphisms. Let $B$ be the subgroup of $G$ consisting of all automorphisms that take $e_i$ to an $R$-linear combination of $e_1,\ldots,e_i$ and $x_i$ to a $K$-linear combination of $x_0,\ldots,x_i$. It is a Borel group of $G$. A submodule $M$ of $F$ is called {\it Borel-fixed}  if $\gamma(M) = M$ for all $\gamma \in B$. A combinatorial characterization of such submodules is given in \cite{Pardue-th}, Proposition II.5. They occur naturally because of the following result. For every graded submodule $M$ of $F$ there is a Zariski open subset $U \subset G$ and a Borel-fixed module $N$ such that $N = \iin (\gamma(M))$ for all $\gamma \in U$ (\cite{Pardue-th}, Example I.7, or modify the proofs in \cite{Eisenbud-book} or \cite{Green-gin} for the case $F = R$). The module $N$ is called the {\it generic initial module} of $M$ and is denoted by $\gin (M)$. 

Generic initial modules can be used to compute the smallest extended degree. This result extends Theorem 3.4.2 in \cite{Gunst}. 

\begin{lemma} \label{lem-gin} 
For every graded submodule $M$ of $F$ we have 
$$
\bdeg (F/M) = \bdeg(F/\gin(M)). 
$$
\end{lemma} 

\begin{proof} 
The proof is analogous to the one of Gunston for the case $F = R$. We sketch it for the sake of completeness. 

We induct on $\dim F/M$. The $0$-dimensional case being clear we may assume $\dim F/M \geq 1$. 
First, suppose that $\depth F/M > 0$. Let $l \in R$ be a general linear form. Then we get by induction 
\begin{eqnarray*} 
\bdeg (F/M) & = & \bdeg (F/M + l F) = \bdeg (F/\gin (M + l F)) \\
& = & \bdeg (F/\gin(M) + x_n F) = \bdeg(F/ \gin(M) + l F) \\
& = & \bdeg(F/\gin(M)) 
\end{eqnarray*} 
where we have also used properties  of generic initial modules under hyperplane sections that are straightforward extensions of Green's results \cite{Green-gin} for generic initial ideals. 

Second, assume $\depth F/M = 0$. We will use the saturation of $M$, i.e. \ 
$$
M^{sat} := \bigcup_{k \geq 1} (M :_F \fm^k). 
$$
Since $\depth F/M^{sat} > 0$ and saturation commutes with formation of generic initial modules 
we obtain 
\begin{eqnarray*} 
\bdeg (F/M) & = & \bdeg(F/M^{sat}) + \len (M^{sat}/M) \\
& = & \bdeg (F/\gin(M)^{sat}) + \len (\gin(M)^{sat}/\gin(M)) \\
& = & \bdeg (F/\gin(M)) 
\end{eqnarray*} 
and we are done. 
\end{proof}

We are ready for the main results of this section. The first one shows that the converse of Theorem \ref{thm-extr-algebras} is true if we use the smallest extended degree. It also gives a nice tie with the theory of Gr\"obner bases. 

\begin{theorem} \label{thm-char-s}
Let $I$ be a homogeneous ideal in $R = K[x_0,\ldots,x_n]$ where $K$ is an infinite field.  Assume that $I$ does not contain any linear form. 
Then the following conditions are equivalent: 
\begin{itemize} 
\item[(a)] $\reg (R/I) = \bdeg (R/I) - 1$ \; and \; $\depth R/I = t$. 
\item[(b)] $I$ is an ideal as in Lemma \ref{lem-ex-ideal} with $n+1-t$ minimal generators. 
\item[(c)] There are integers $d_t,\ldots,d_n \geq 0$ with $0 \leq t \leq n$, $d_n, d_t >0$ such that 
$$
\gin (I) = (x_0^{d_n+1}, x_0^{d_n} x_1^{d_{n-1} + 1},\ldots, x_0^{d_n} \dots x_{n-t-2}^{d_{t+2}} x_{n-t-1}^{d_{t+1}+1}, x_0^{d_n} \dots x_{n-t-1}^{d_{t+1}} x_{n-t}^{d_{t}}). 
$$ 
\end{itemize}
\end{theorem} 

\begin{proof} 
Theorem \ref{thm-extr-algebras} shows that (a) implies (b). In order to prove the converse we have to show for every ideal $I$ as in (b) that 
$$
\bdeg (R/I) = d_t + \ldots + d_n 
$$
because $\reg (R/I) = d_t + \ldots + d_n - 1$ as a consequence of Lemma \ref{lem-ex-ideal}. 

We will induct on $n$. Condition (1) for $\bdeg$ and Lemma \ref{lem-ex-ideal} provide 
\begin{eqnarray*}
\bdeg (R/I) & = & \bdeg (R/I^{sat}) + \len (\HH^0 (R/I)) \\
& = & \bdeg (R/I^{sat}) + d_0. 
\end{eqnarray*} 
In particular, this implies the claim for $n=0$, thus we may assume $n \geq 1$. 
Let $s = \min \{i \geq 1 \s d_i \geq 1\}$. This minimum exists because $d_n \geq 1$. Then, using Lemma \ref{lem-ex-ideal} we see that 
$$
I^{sat} = (f_n l_0, f_n f_{n-1} l_1,\ldots, f_n \ldots f_{s+1} l_{n-s-1}, f_n \ldots f_s)
$$
is again an ideal as in Lemma \ref{lem-ex-ideal} where $\depth R/I^{sat} = s$. Hence, it remains to consider the case where $t = \depth R/I > 0$. Changing coordinates we may assume that $x_n$ is general for $R/I$ and we have 
$$
\bdeg (R/I) = \bdeg (R/I + x_n R). 
$$ 
But $R/I + x_n R \cong \bar{R}/\bar{I}$ with $\bar{R} = R/x_n R$ and 
$$
\bar{I} = (g_{n-1} \bar{l}_0, g_{n-1} g_{n-2} \bar{l}_1,\ldots, g_{n-1} \ldots g_{t} \bar{l}_{n-t-1}, g_{n-1} \ldots g_{t-1})
$$ 
where $g_i$ is the residue class of $f_{i+1}$ and $\bar{l_i}$ is the residue class of $l_i$. We may apply the induction hypothesis to $\bar{I}$. Hence we obtain 
$$
\bdeg (\bar{R}/\bar{I}) = d_t + \ldots + d_n 
$$ 
which provides  $\bdeg (R/I) = d_t + \ldots + d_n $, as claimed. 
\smallskip 

Now, we show that (c) implies (a). Note that the given $\gin (I)$ is an ideal as in Lemma \ref{lem-ex-ideal}. Hence, we have just seen that $t = \depth R/\gin(I)$ and 
$$
\reg(R/\gin(I)) = \bdeg(R/\gin(I)) -  1. 
$$
According to \cite{BS-invent} we have 
$$
\reg (R/I) = \reg(R/\gin(I)). 
$$ 
Combined with Lemma \ref{lem-gin}  and $\depth R/I = \depth R/\gin(I)$ we get the equations in (a). 

Reversing the last argument we see that (a) implies the fact that $\gin(I)$ is an ideal as in Lemma \ref{lem-ex-ideal} and $\depth R/\gin(I) = t$. Since $\gin(I)$ is a Borel-fixed ideal it must be of the form as claimed. 
\end{proof} 

\begin{remark} \label{rem-embedded}
Unless the ideal is principal the ideals occurring in Theorem \ref{thm-char-s}(ii) are not unmixed. However, the lower dimensional components are not necessarily embedded.  
\end{remark} 

Combined with Lemma \ref{lem-ex-ideal}, Theorem \ref{thm-char-s} says precisely for which Hilbert series $H$ where is a graded $K$-algebra $A$ with $H_A = H$ and $\reg (A) = \bdeg (A) - 1$. A lot more is true. There is another purely numerical characterization of these algebras. 

\begin{theorem} \label{thm-h-func} 
Let $A \neq K$ be a standard graded $K$-algebra where $K$ is an infinite field.  Then the following conditions are equivalent 
\begin{itemize} 
\item[(a)] $A$ has maximal regularity with respect to the smallest extended degree, i.e.\  
$$
\reg (A) = \bdeg (A) - 1. 
$$
\item[(b)]  The Hilbert series of $A$ is  
$$
H_A (z) = \sum_{j=0}^n \frac{1-z^{d_j}}{(1-z)^{j+1}} \prod_{i=j+1}^n z^{d_i} 
$$ 
for some integers $d_0,\ldots,d_n \geq 0 $  where $n := \dim A$ and  $d_n > 0$.  
\end{itemize} 
\end{theorem} 

\begin{remark} \label{rem-num-char} 
Using Theorem \ref{thm-char-s} we see that the integers $d_i$ occurring in condition (b) have a cohomological interpretation. Indeed, Lemma \ref{lem-ex-ideal} shows that $\HH^i(A)^{\vee}$ is (up 
to a degree shift) a hypersurface ring of degree $d_i$. Here degree zero means that the cohomology module is trivial. 
\end{remark} 

\begin{proof}[Proof of Theorem \ref{thm-h-func}] 
One direction is clear. Theorem \ref{thm-char-s} in conjunction with Lemma \ref{lem-ex-ideal} shows that (a) implies (b). 

For the other implication we prove by induction on $n$ that (b) has as a consequence 
$$
\reg (A) = \bdeg (A) - 1 = d_0 + \ldots + d_n - 1. 
$$ 
This is clear if $n = 0$. Let $n \geq 1$. Notice that the Hilbert series of $A$ conveys that $d_n = \deg A$ and $n+1 = h_A (1) = \rank [A]_1$. Thus, we see that $A \cong R/I$ where $R = K[x_0,\ldots,x_n]$. Since $\dim A = n$ the ideal $I$ has codimension one in $R$. Hence, the greatest common divisor of the elements in $I$ is a homogeneous polynomial $f$ of degree $d_n = \deg A$. Setting $B := R/I : f$, multiplication by $f$ on $A$ provides the exact sequence 
\begin{equation} \label{eq-seqq} 
0 \to B (-d_n) \to A \to R/f R \to 0. 
\end{equation} 
Hence the Hilbert series of $B$ is 
$$
H_B (z) = \sum_{j=0}^{n-1} \frac{1-z^{d_j}}{(1-z)^{j+1}} \prod_{i=j+1}^{n-1} z^{d_i} 
$$
and the induction hypothesis provides 
$$
\reg (B) = \bdeg (B) - 1 = d_0 + \ldots + d_{n-1} - 1. 
$$
Since $\reg (R/f R) = d_n -1$ we obtain from the sequence (\ref{eq-seqq}) 
$$
\reg (A) = \reg (B) + d_n = d_0 + \ldots + d_n - 1.
$$
It remains to compute $\bdeg (A)$. Due to Theorem \ref{thm-char-s}, Lemma \ref{lem-ex-ideal} applies to $B$ and shows that the Hilbert series of $\HH^0(B)$ is 
$$
H_{\HH^0(B)} (z) = \frac{1-z^{d_0}}{1-z} \prod_{i=1}^n z^{d_i}. 
$$ 
Using the sequence (\ref{eq-seqq}) again we get that $\HH^0(A) \cong \HH^0(B)$. Therefore the Hilbert series of $C := A/\HH^0(A)$ is 
\begin{eqnarray*} 
H_C (z) & = & H_A (z) - H_{\HH^0(B)} (z) \\
& = & \sum_{j=1}^n \frac{1-z^{d_j}}{(1-z)^{j+1}} \prod_{i=j+1}^n z^{d_i}. 
\end{eqnarray*} 
Let $l \in R$ be a general linear form. It is not a zero divisor on $C$. Hence we get for the Hilbert series of $C/l C$ 
$$
H_{C/l C} (z) = (1 - z) \cdot H_C (z) = \sum_{j=1}^{n} \frac{1-z^{d_j}}{(1-z)^{j}} \prod_{i=j+1}^n z^{d_i}. 
$$ 
Thus, we can apply the induction hypothesis to $C/lC$ and obtain
$$
\bdeg (C) = \bdeg (C/ l C) = d_1 + \ldots + d_n. 
$$
Since $\HH^0(A)$ has length $d_0$ it follows 
$$
\bdeg (A) = \bdeg (C) + \len (\HH^0(A)) = d_0 + \ldots + d_n 
$$ 
completing the proof. 
\end{proof} 

Expressing the information on the Hilbert series of the algebras more explicitly by using the Hilbert function illustrates again how special these algebras are. We work this out in low dimensions. 

\begin{example} \label{ex-h-func}  Let $A$ be an algebra as in Theorem \ref{thm-h-func} (b). 

(i) If $A$ is artinian then the Hilbert function of $A$ is $1,1,\ldots,1,0,0 \ldots$. 

(ii) If $\dim A = 1$ then the Hilbert function is 
$$ 
h_A (j) = \left \{ 
\begin{array}{cl} 
j + 1 & \mif 0 \leq j \leq d_1 - 1 \\[1ex]
d_1 + 1 & \mif d_1 \leq j \leq d_0 + d_1 -1 \\[1ex]
d_1 & \mif d_0 + d_1 \leq j. 
\end{array} \right.
$$
\end{example} 

Using Lemma \ref{lem-extr-is-cyclic} it is possible to  extend the
previous results to modules. For example,  the  module version of
Theorem \ref{thm-h-func} is.  

\begin{corollary} \label{cor-h-mod} Let $R$ be a noetherian 
  polynomial ring over the infinite field $K$, and
  let $M \neq 0$ be a finitely generated graded $R$-module. Then the
  following conditions are equivalent  
\begin{itemize} 
\item[(a)] $M$ has maximal regularity with respect to the smallest
  extended degree, i.e.\   
$$ 
\reg (M) = e^+_R(M) + \bdeg (M) - 1. 
$$
\item[(b)] $M$ is a cyclic $R$-module with  Hilbert series  
$$
H_M (z) = z^k \cdot \left [\sum_{j=0}^n \frac{1-z^{d_j}}{(1-z)^{j+1}}
  \prod_{i=j+1}^n z^{d_i} \right ] 
$$ 
for some integers $k, d_0,\ldots,d_{n}$ where $n = \dim M$, $d_0,\ldots,d_{n-1} \geq 0$,  and  $d_n > 0$.  
\end{itemize} 
\end{corollary}  

\begin{proof} 
Using $H_{M(-k)}(z) = z^k \cdot H_M(z)$ this follows by combining Lemma \ref{lem-extr-is-cyclic} and Theorem \ref{thm-h-func}. 
\end{proof} 


\section{Comments and open problems} \label{sec-open} 

We begin by discussing the assumption about the ground field that we have made  in the previous sections. Then we show that our results broadly generalize the main result of St\"uckrad and Vogel in \cite{St-V-mult}.

\begin{remark} \label{rem-finite-field} 
The assumption  that the ground field $K$ is infinite was necessary  in order to guarantee that there are filter regular elements. However, typically one would expect
  that a degree does not change when the base field $K$ is
  extended. If we add this requirement as an additional condition for an extended degree then Theorem \ref{thm-reg-bound} and Corollary \ref{cor-reg-b} become valid over arbitrary ground fields. Indeed, we could replace $K$ by its algebraic closure if necessary. Note that Vasconcelos' homological degree satisfies the extra requirement above. 
\end{remark} 
\smallskip 

Every extended degree $\Deg$  satisfies $\Deg (M) \geq \deg (M)$ and quality holds if and only if $M$ is a Cohen-Macaulay module. Following \cite{DG-Vasc}, the difference 
$$
I(M) := \Deg (M) - \deg (M) 
$$
is called the {\it Cohen-Macaulay deviation} of $M$. Using it we can  state another consequence of Theorem \ref{thm-hypsec}. 

\begin{proposition} \label{prop-buchs} 
Let $M$ be a finitely generated graded $R$-module of dimension $d \geq 1$. Let $\fq \subset R$ be a parameter ideal generated by $d$ sufficiently general linear forms. Suppose that the ground field $K$ is infinite. Then we have for any extended degree
\begin{equation*} \label{eq-buchs}
\len (M/ \fq M) - e_0 (\fq; M) - \reg (M/\fq M) \leq I(M) - \reg (M) 
\end{equation*} 
where $e_0(\fq; M)$ denotes the multiplicity of $\fq$ (cf., e.g., \cite{ZS}, Chapter VIII). 
\end{proposition} 

\begin{proof} 
Successive application of Theorem \ref{thm-hypsec} provides 
\begin{eqnarray*} 
\len (M/\fq M) - \reg (M/\fq M) & \leq & \Deg (M) - \reg (M) \\
& = & \deg (M) + I(M) - \reg (M). 
\end{eqnarray*} 
Since the minimal generators of $\fq$ are sufficiently general by assumption, they form an $M$-filter regular sequence. Thus, \cite{Ausl-B}, Proposition 4.7, shows that they form a reducing system of parameters for $M$. Hence, we have $e_0(\fq; M) = \deg (M)$. Combined with the estimate above the claim follows. 
\end{proof} 

It is useful to specify what means ``sufficiently general'' for the linear forms in the last result provided we deal with a specific extended degree. All we need is to know which assumptions on the linear form $l$ guarantee that $\Deg (M/l M) \leq \Deg (M)$. This follows by an analysis of the proof of Theorem \ref{thm-hypsec}.

In case of the homological degree this has been achieved by Vasconcelos. Following \cite{Vasc-hom-deg}, Definition 2.12, a linear form $l \in R$ is called a {\it special hyperplane section} of $M$ if it is filter regular for all the iterated extension modules 
$$
\Ext_R^{i_1}(\Ext_R^{i_2}(\ldots(\Ext_R^{i_p}(M, R), R),\ldots, R) 
$$
where $i_1 \geq i_2 \geq \ldots \geq i_p \geq 0$ are any integers. 

Suppose that the module $M$ has dimension $d$. We call the ideal $\fq = (l_1,\ldots,l_d) \subset R$ a {\it special parameter ideal} of $M$ if  $l_i$ is a special hyperplane section of $M/(l_1,\ldots,l_{i-1}) M$ for all $i = 1,\ldots,d$. 

In case of the homological degree we denote the Cohen-Macaulay deviation of the module $M$ by 
$$
I_h (M) = \hdeg (M) - \deg (M). 
$$ 
Using this notation Proposition \ref{prop-buchs} provides: 

\begin{corollary} \label{cor-bu}
Let $\fq \subset R$ be a special parameter ideal of the module
$M$. Then we have for the homological degree  
\begin{equation*} \label{eq-buchs-cor} 
\len (M/ \fq M) - e_0 (\fq; M) - \reg (M/\fq M) \leq I_h (M) - \reg (M). 
\end{equation*} 
\end{corollary} 

\begin{proof} 
Vasconcelos' Theorem 2.13 in \cite{Vasc-hom-deg} shows for a special
hyperplane section $l$ of $M$  
$$
\hdeg (M/l M) \leq \hdeg (M). 
$$ 
As pointed out above, the claim then follows by Proposition \ref{prop-buchs}. 
\end{proof} 

The last result can be further simplified if $M$ is a module with
finite cohomology, i.e.\ if all the modules $\HH^i(M)$, $0 \leq i <
\dim M$, have finite length. Note that this is the case if and only if
$M$ is an equidimensional locally Cohen-Macaulay module. It is easy to
see that then a linear form is a special hyperplane section of the
module $M$ if and only if it is a parameter for $M$. Thus, the last
corollary implies:  

\begin{corollary} \label{cor-fcoh} 
Let $M$ be a an equidimensional locally Cohen-Macaulay graded
module. Let $\fq \subset R$ be an ideal generated by a system of
linear parameters of $M$ and suppose that $K$ is infinite. Then we
have  
\begin{equation*} 
\len (M/ \fq M) - e_0 (\fq; M) - \reg (M/\fq M) \leq I_h (M) - \reg (M). 
\end{equation*}
\end{corollary} 

Note that with the assumptions of the last result, the Cohen-Macaulay
deviation is  
$$
I_h (M) = \sum_{0}^{d-1} \binom{d-1}{i} \len (\HH^i(M))  
$$
where $d = \dim M$. Hence, Corollary \ref{cor-fcoh} is a
generalization of the Main theorem in \cite{St-V-mult} (cf.\ also
\cite{Ach-Sch}, \cite{N-Nachr}) from algebras to modules. There
St\"uckrad and Vogel have also considered the problem for which
algebras we have equality in Corollary \ref{cor-fcoh} for every ideal
$\fq$ under consideration. Using Corollary \ref{cor-equal-in-thm}
their methods can be extended to the case of modules and provide, for
example, the following result.  If $M$ has positive depth then we have
equality in Corollary \ref{cor-fcoh} for every ideal $\fq$ under
consideration if and only if $M$ is a Buchsbaum module.  

It seems interesting to analyze the analogous problem in the more
general situation of Corollary \ref{cor-bu}. We pose this as an open
problem.


\end{document}